\documentclass[11pt]{article}
\usepackage{amsmath}
\newcommand{\integ}[1]{\ensuremath{\int_{0}^{T}\!\!\int_{\Omega}{#1}
    \,dx\,dt}} 
\newcommand{\inte}{\ensuremath{\int_{0}^{T}\!\!\int_{\Omega}}} 
\newcommand{\Ma}[3]{\ensuremath{|\nabla {#1}|^{#3}\nabla{#1} -
    |\nabla{#2}|^{#3} \nabla{#2}}} 
\newtheorem{thm}{Theorem}
\newtheorem{lemma}[thm]{Lemma}

\DeclareMathOperator{\divergence}{div}
\DeclareMathOperator*{\esssup}{ess\,sup}

\begin{document}
\title{On the Time Derivative in a Quasilinear Equation}
\author{Peter Lindqvist}
\date{Norwegian University of Science and Technology}
\maketitle
\begin{center}
\textsf{Dedicated to the memory of \emph{Juha Heinonen} 1960--2007}
\end{center}
\maketitle
\section{Introduction}
 
The regularity theory for certain parabolic differential equations of
the type 
\begin{equation}\frac{\partial u}{\partial t} =
 \divergence\textbf{A}(x,t,u,\!\nabla u)
\end{equation}
does often not treat the time derivative 
 $u_{t}$, which is regarded as
a distribution. Thus the time derivative is a neglected object.
In this note we will prove that the weak solutions of the Evolutionary
p-Laplace Equation
\begin{equation}\frac{\partial u}{\partial t} = \divergence(|\nabla
  u|^{p-2}\nabla u) 
\end{equation}
have a time derivative $u_{t}$ in Sobolev's sense. In particular, 
$u_{t}$ is not merely a distribution but a measurable function,
belonging to some space $L^{q}_{\rm{loc}}$. 

No doubt, analogous results are known to the experts. The evident fact is
that, if the right-hand side of the equation (the divergence part)
is a function, so is the left-hand side (the time derivative). Indeed,
it has been noted that this yields a derivative even for systems, as in
section 7 of ~\cite{AM}, and frequently the required estimates appear at
intermediate steps in advanced proofs aiming at the continuity of the
gradient $\nabla u$, as in ~\cite{Y}. For equation (2) much simpler proofs
are accessible. It is an advantage to have the time derivative at ones
disposal at an early stage of the theory. Therefore we have found it
worth our while to present a direct and succinct proof of the
existence and summability of the time derivative. We are able to avoid the
use of Moser's and de Giorgi's iterations. Deeper regularity
properties are beyond the scope of this note.

The Evolutionary p-Laplace Equation is \emph{degenerate} for $p >
2$ and \emph{singular} for $1 < p < 2$. We will restrict
ourselves to the cases $2 \leq p < \infty$. We refer to the books ~\cite{dB}
and  ~\cite{WZ} about this equation\footnote{It is also called the
  non-Newtonian equation of filtration.}, originally encountered more than half
a century ago by Barenblatt. The proof can readily be extended to
equations like $$\frac{\partial u}{\partial t} =
\sum_{i,j}\frac{\partial}{\partial x_{i}} \left(\left|\sum_{k,m}
      a_{k,m}\frac{\partial u}{\partial x_{k}}\frac{\partial
        u}{\partial x_{m}}\right|^{\frac{p-2}{2}}\!
    a_{i,j}\frac{\partial u}{\partial x_{j}}\right)$$
provided that the constant matrix $(a_{i,j})$ satisfies the
ellipticity condition $$ \sum a_{i,j}\xi_{i}\xi_{j} \geq \lambda |\xi|^{2}
.$$ Also the case $ a_{i,j} =  a_{i,j}(t)$ is easy, but further
generalizations seem to require more refined assumptions. The result
is not valid for all equations of the type (1)\footnote{This may
  explain why the time derivative is neglected in the literature.}.
Here we are content
with the more pregnant formulation in terms of the Evolutionary
p-Laplace Equation.

\section{The Caccioppoli Estimate}

We first define the concept of solutions, then we state the main
theorem. The rest of the section is devoted to a Caccioppoli
estimate.

Suppose that $\Omega$ is a bounded domain in $\mathbf{R^{n}}$ and
consider the space-time cylinder $\Omega_{T} = \Omega \times
(0,T)$. In the case $p \geq 2$ we say that $u \in
L^{p}(0,T;W^{1,p}(\Omega))$ is a \emph{weak solution} of the
Evolutionary p-Laplace Equation, if 
\begin{equation}
\integ{(-u\phi_{t} + \langle |\nabla u|^{p-2},\!\nabla \phi \rangle)\,} = 0
\end{equation}
for all $ \phi \in C^{1}_{0}(\Omega_{T})$. (The singular case $1 < p <
2$ requires an extra \emph{a priori} assumption, for example, $u \in
L^{\infty}(0,T;L^{2}(\Omega))$ will do.) In particular, one has
$$ \integ{(|u|^{p} + |\nabla u|^{p})} <  \infty.$$
By the regularity theory one may regard $u(x,t)$ as continuous, a fact
which we need not use. The main result is the following.
\begin{thm}
Let $2\leq p <\infty$. If $u = u(x,t)$ is a weak solution, then the
time derivative $u_{t}$ exists (in Sobolev's sense) and
  $u_{t} \in L^{p/(p-1)}_{\rm{loc}}(\Omega_{T})$.
\end{thm}

The proof is based on the applicability of the rule 
\begin{equation}
\integ{u\phi_{t}} = - \integ{\phi\,\nabla \! \cdot \!(|\nabla u|^{p-2}\nabla
  u)}
\end{equation}
when $\phi \in C^{1}_{0}(\Omega_{T})$. Thus the theorem follows
provided that it first be properly established that the Sobolev
derivatives $\partial/\partial x_{j}(|\nabla u|^{p-2}\nabla u)$, 
appearing in the formula, exist
and belong to $L^{p/(p-1)}_{\rm{loc}}(\Omega_{T})$. The main task is thus
to prove differentiability in \emph{the $x$-variable}.

To begin with, we need a variant of the Caccioppoli estimate for the
difference $u(x+h,t) - u(x,t)$ where $h$ is a small increment in the
desired direction. If $\phi$ is a given test function with compact
support, then also the translated function $v(x,t) = u(x+h,t)$ is a
weak solution in some subdomain containing the support of $\phi$,
provided that $|h|$ is small enough. Subtracting the equations for
$u(x,t)$ and $u(x+h,t)$ we obtain  

\begin{eqnarray}
\integ{\langle\Ma{u(x\!+\!h,t)}{u(x,t)}{p-2},\nabla \phi(x,t)
  \rangle}\nonumber \\
= \integ{(u(x+h,t)-u(x,t))\phi_{t}(x,t)}.
\end{eqnarray}
Choose the test function $$\phi(x,t) =
\eta(t)\zeta(x)^{p}(u(x+h,t)-u(x,t))$$ where $\zeta \in
C^{\infty}_{0}(\Omega)$, $0\leq \zeta(x) \leq 1$, and $\eta(t)$ is a
cut-off function, $0 \leq \eta(t) \leq 1$, and $\eta(0) = \eta(T) =
0.$ Strictly speaking it is not an admissible one, because $\phi_{t}$
contains the forbidden time derivative $u_{t}$. A \emph{formal}
calculation yields the \emph{Caccioppoli estimate}\footnote{This is a
  slight abuse of the name, since there is no estimate yet.}
\begin{multline*}
\inte \eta(t)\zeta(x)^{p}\langle \Ma{u(x\!+\!h,t)}{u(x,t)}{p-2},\\
\shoveright{\nabla u(x\!+\!h,t) - \nabla u(x,t) \rangle\,dx\,dt}\\
= - p \inte\eta(t)\zeta(x)^{p-1}(u(x\!+\!h,t)-u(x,t))\\
\shoveright{\times \langle\Ma{u(x+h,t)}{u(x,t)}{p-2},
\nabla \zeta(x,t)\rangle \,dx\,dt}\\
+ \frac{1}{2} \integ{\eta'(t)\zeta(x)^{p}(u(x\!+\!h,t)-u(x,t))^{2}}
\end{multline*}
after some integrations by part of the integral containing $\phi_{t}.$

In order to justify the use of the test function above we introduce
the convolution $$
f(x,t)^{*} = \inte{f(x-y,t-\tau)\rho_{\sigma}(y,\tau)\,dy\,d\tau},$$
where $\rho_{\sigma}$ is a smooth non-negative function with compact
support in the ball $|y|^{2} + \tau^{2} \leq \sigma^{2}$;
$\sigma$ is small. (In fact, convolution only in the time variable
would suffice. The familiar Steklov average works well.) With the
abbreviations
 $u = u(x,t)$ and $v = u(x+h,t)$
we obtain the averaged identity
\begin{eqnarray}
\integ{\langle(|\nabla{v}|^{p-2}\nabla{v})^{*}
  -(|\nabla{u}|^{p-2}\nabla{u})^{*},\nabla \phi \rangle} \nonumber\\
= \integ{(v^{*} - u^{*}) \phi_{t}}
\end{eqnarray}
from equation (5). This is a standard procedure. The parameter
$\sigma$ has to be less than a bound depending on $|h|$ and on the
distance from the support of the test function $\phi$ to the
boundary. Now we insert the test function
$$ \phi(x,t) = \eta(t)\zeta(x)^{p}(v(x,t)^{*}-u(x,t)^{*})$$
into (6). This is an admissible one. Again the integral containing
$\phi_{t}$ becomes $$
\frac{1}{2} \integ{\eta'(t)\zeta(x)^{p}(v(x,t)^{*}-u(x,t)^{*})^{2}}.$$
Here we may safely let $\sigma \rightarrow 0$. The terms coming from
$\nabla \phi$ cause no problem, when $\sigma \rightarrow 0$. Thus we
arrive at the Caccioppoli estimate again, but this time the procedure
was duly justified.

\section{Estimation of Difference Quotients}

We aim at proving differentiability in the variable $x$ of the
auxiliary vector field
$$F(x,t) = |\nabla u(x,t)|^{(p-2)/2}\nabla u(x,t)$$
by bounding its integrated difference quotients. Notice that we have
$(p-2)/2$ in place of the desired exponent $p-2$, the transition to
which is explained in section 4. In the stationary case this expedient
quantity was employed by Bojarski and Iwaniec, cf. ~\cite{BI}. They used the
elementary inequalities
\begin{equation} 
\frac{4}{p^{2}}\left||b|^{\frac{p-2}{2}}b -|a|^{\frac{p-2}{2}}a
\right|^{2} \leq \\\langle |b|^{p-2}b - |a|^{p-2}a, b-a \rangle
\end{equation}
\begin{equation}
\left||b|^{p-2}b -|a|^{p-2}a\right| \leq \\
(p-1) \left(|b|^{\frac{p-2}{2}} + |a|^{\frac{p-2}{2}}\right)
\left||b|^{\frac{p-2}{2}}b -|a|^{\frac{p-2}{2}}a\right|
\end{equation}
for vectors, where $p\geq 2.$\footnote{A proof is worked out in ~\cite{L}.}

The partial differentability of $F$ often comes as a by--product of
more advanced considerations aiming at establishing the continuity of 
$\nabla u$ itself, as, for example, in ~\cite{Y}. We give a simpler proof
below, avoiding iterations. (Needless to say, we do not reach the
continuity of $\nabla u$ this way.) We write $DF$ for the 
matrix with the elements $$\frac{\partial}{\partial x_{j}}\!\left(
 |\nabla u|^{(p-2)/2}\frac{\partial u}{\partial x_{i}}\right).$$
\begin{lemma}
Let $p > 2$. The derivatives $DF$ exist in Sobolev's sense and $DF
\in L^{2}_{loc}(\Omega_{T})$. The estimate 
\begin{eqnarray}
\int_{\tau}^{T}\!\!\int_{\Omega}\zeta(x)^{p}|DF|^{2}\,dx\,dt \leq
\frac{c}{\tau}\int_{0}^{\tau}\!\!\int_{\Omega}\zeta(x)^{p}|\nabla
u(x,t)|^{2}\,dx\,dt
\nonumber \\
+ \integ{(\zeta(x)^{p}+|\nabla \zeta(x,t)|^{p})|\nabla u(x,t)|^{p}}
\end{eqnarray}
holds when $\tau > 0$. Here $\zeta \in C_{0}^{\infty}(\Omega)$,
 $\zeta(x) \geq 0$.
\end{lemma}

\emph{Proof:} Proceeding from the Caccioppoli estimate in section 2 we obtain,
using the elementary inequalities (7) and (8), 
\begin{multline}
\frac{4}{p^{2}} \integ{\eta(t) \zeta(x)^{p}|F(x+h,t)-F(x,t)|^{2}}
\\
\leq  \frac{1}{2} \integ{\eta'(t) \zeta(x)^{p}
  (u(x+h,t)-u(x,t))^{2}}\\
+ p(p-1) \inte
\left(\eta(t)^{\frac{1}{2}}\zeta(x)^{\frac{p}{2}}|F(x\!+\!h,t)-F(x,t)|\right)
\\ \times \left(\eta(t)^{\frac{1}{2}}|u(x\!+\!h,t)-u(x,t)||\nabla \zeta(x)|
\right)\phantom{ab} \\
\times \left(\left|\nabla u(x+h,t)\right|^{\frac{p-2}{2}}+
  \left|\nabla u(x,t)\right|^{\frac{p-2}{2}}
\right)\zeta(x)^{\frac{p-2}{2}} \,dx\,dt. 
\end{multline}
Divide both sides by $|h|^{2}$ and use the inequality
$$abc \leq \frac{\varepsilon^{2}a^{2}}{2}
+\frac{\varepsilon^{-p}\,b^{p}}{p}+  \frac{(p-2)c^{2p/(p-2)}}{2p} $$
where the exponents $2, p, 2p/(p-2)$ are conjugated. It follows that
the last integral is majorized by
\begin{eqnarray*}
\frac{p(p-1)\varepsilon^{2}}{2}\integ{\eta(t) \zeta(x)^{p}
\left|\frac{F(x+h,t)-F(x,t)}{h}\right|^{2}}\\
+ (p-1)\varepsilon^{-p}\integ{\eta(t)^{\frac{p}{2}}
  \left|\frac{u(x\!+\!h,t)-u(x,t)}{h}
\right|^{p}|\nabla \zeta(x,t)|^{p}}\\
+ c_{p} \integ{\zeta(x)^{p}(|\nabla u(x+h,t)|^{p} +|\nabla u(x,t)|^{p})}\rlap{.}
\end{eqnarray*}
Choose $\varepsilon > 0$ so small that the term with
 $\varepsilon^{2}$ is absorbed by the left--hand side of (10), for
example, take $p(p-1)\varepsilon^{2}/2 = 2/p^{2}$, which is half of $4/p^{2}$.Then
\begin{eqnarray*}
\frac{2}{p^{2}}\integ{\eta(t) \zeta(x)^{p}\left|
\frac{F(x+h,t)-F(x,t)}{h}\right|^{2}}\\
\leq  \frac{1}{2} \integ{\eta'(t) \zeta(x)^{p}
  \left|\frac{u(x+h,t)-u(x,t)}{h}\right|^{2}}\\
+ a_{p}\integ{\left|\frac{u(x\!+\!h,t)-u(x,t)}{h}
\right|^{p}|\nabla \zeta(x,t)|^{p}}\\
+ c_{p} \integ{\zeta(x)^{p}(|\nabla u(x+h,t)|^{p} +|\nabla u(x,t)|^{p})}\rlap{.}
\end{eqnarray*}
Let us finally select $\eta(t)$ as a piecewise linear cut-off function
so that $\eta(t) = 1$ when $\tau \leq t \leq T-\beta$. Since 
$\eta'(t) < 0$, when $t > T - \beta$, we may omit that portion of the
integral in question and then let $\beta \rightarrow 0$. There is no
trace left of $\beta$ in the formula. We can further arrange it so
that $\zeta(x) = 1$ in an arbitrary compact subset of $\Omega$. The
characterization of Sobolev's space in terms of integrated difference
quotients guarantees that the derivatives $DF$ exist. As $h
\rightarrow 0$ we arrive at the desired estimate. This concludes the
proof.

\medskip\noindent
\emph{Remark:} Using the 'lost' interval $[T-\beta,T]$ in an effective
way, a standard procedure yields an estimate also of
$$\esssup_{0 < t < T}\int_{\Omega}\zeta(x)^{p}|\nabla
u(x,t)|^{2}\,dx.$$
\section{The end of the proof in the degenerate case and comments on
  the singular case}

We are in the position to conclude the proof in the case $p > 2$. We
have \begin{eqnarray*}|F|^{2} = |\nabla u|^{p},& |\nabla u|^{p-2}\nabla u =
|F|^{1-2/p}F.\end{eqnarray*}
Thus, in virtue of the lemma,
\begin{equation}
\left|\frac{\partial}{\partial x_{j}}(|\nabla u|^{p-2}\nabla u)\right|
\leq 2\left(1-
\frac{1}{p}\right) \left|F \right|^{\frac{p-2}{p}} \left|\frac{\partial
  F}{\partial x_{j}}\right|
\end{equation}
and, by H\"{o}lder's inequality,
$$\frac{\partial}{\partial x_{j}}(|\nabla u|^{p-2}\nabla u)\, \in \,
L^{\frac{p}{p-1}}_{\rm{loc}}
(\Omega_{T})$$
because $F \in L^{2}(\Omega_{T})$ and $DF \in
L^{2}_{loc}(\Omega_{T})$. Finally, the theorem follows from the rule
(4). This concludes the proof.

\medskip\noindent
\emph{Remark:} In fact, $F \in L^{\infty}_{loc}(\Omega_{T})$ and hence
one can prove that $u_{t} \in L^{2}_{\rm{loc}}(\Omega_{T})$, which is
stronger. However, this boundedness of $F$ requires more advanced
regularity theory. For example, in ~\cite{dB} the continuity of $\nabla u$,
and consequently of $F$, is proved. 

\medskip
Let us finally mention that in the \emph{singular case} $1 < p <2$ one
rather easily obtains that the Sobolev derivatives $u_{x_{i}x_{j}}$
of the second order and $DF$ exist and belong to
$L^{2}_{loc}(\Omega_{T})$. (When $p > 2$, $u_{x_{i}x_{j}}$is more difficult to
achieve!) Unfortunately, one encounters a new complication in (11),
caused by the negative exponents. Thus the full regularity theory
seems to be needed. In section 2 of ~\cite{Y} the crucial estimate 
$$\iint |\nabla
u|^{2(p-2)}|D^{2}u|^{2} \,dx\,dt < \infty$$
is given for the range $2 \geq p > \max[3/2,2n/(n+2)]$. To this one
may add that the range $1 < p < 2n/(n+2)$ is not well understood in
general.

\end{document}